\documentclass[12pt,
]{article}
\usepackage[T2A]{fontenc}
\usepackage[cp1251]{inputenc}
\usepackage[english,russian]{babel}

\usepackage{amsfonts,amsmath,amssymb,mathrsfs,amscd}
\usepackage{amsthm}
\usepackage{verbatim}
\usepackage{eucal}
\usepackage{graphicx}
\usepackage{enumerate}
\usepackage{xcolor}

\usepackage[pdftex,unicode,colorlinks,linkcolor=blue,citecolor=red,bookmarksopen,
pdfhighlight=/N]{hyperref}

\theoremstyle{plain}

\newtheorem*{ThL1}{Теорема Лиувилля для целых функций конечного порядка}
\newtheorem*{ThL}{Теорема Лиувилля}

\theoremstyle{definition}
\newtheorem*{remark}{Замечание}

\newcommand{\rad}{\text{\tiny\rm rad}}

\def\RR{\mathbb R}
\def\CC{\mathbb C}

\DeclareMathOperator{\dd}{\,{\rm d\!}}

\begin{document}

\title{К теореме Лиувилля\\ для  целых функций конечного порядка}

\author{Б.\,Н.~Хабибуллин}

\date{29.08.2020}
\maketitle

Основа заметки ---  классическая для {\it целых,\/} или голоморфных  на {\it комплексной плоскости\/} $\CC$, функций    
\begin{ThL} Ограниченная  целая  функция  на  $\mathbb C$ постоянна.   
\end{ThL}

Недавно в работе \cite[лемма 4.2]{BarBelBor18} была дана  версия теоремы Лиувилля для  целых функций {\it конечного порядка\/} на $\CC$, удовлетворяющих по определению ограничению
\begin{equation}\label{fo}
\limsup_{z\to \infty}\frac{\ln\Bigl(1+ \max\bigl\{\ln |f|(z),0\bigr\}\Bigr)}{\ln |z|}<+\infty. 
\end{equation}

 В  \cite[лемма 4.2]{BarBelBor18a} её доказательство откорректировано, а перед её формулировкой в  \cite[преамбула теоремы 2.1]{AlemBarBelHed20} отмечается, что установлена она А.\,А.~Боричевым. 

\begin{ThL1}[{\cite[лемма 4.2]{BarBelBor18}, \cite[лемма 4.2]{BarBelBor18a}, \cite[теорема 2.1]{AlemBarBelHed20}}]  Если целая функция  $f$  конечного порядка  ограничена вне  
множества $E\subset \CC$ нулевой плоской плотности по плоской мере Лебега $\lambda$ в том смысле, что определён предел
\begin{equation}\label{m00}
\lim_{r\to +\infty} \frac{\lambda \bigl(\{z\in E\colon |z|\leq r\}\bigr)}{r^2}=0, 
\end{equation}
то $f$ --- постоянная функция.
\end{ThL1}
Приведённые в  \cite{BarBelBor18} и \cite{BarBelBor18a} её доказательства используют достаточно серьёзный математический аппарат.  Дадим  ниже её  краткое и независимое 
\begin{proof} Рассмотрим субгармоническую  функцию $u:=\ln|f|$ и субгармоническую положительную функцию $u^+:=\sup\{u,0\}$ на $\CC$ конечного порядка в том смысле, что  после замены 
$\ln|f|$ на $u$ и $u^+$ верхний предел в \eqref{fo} по-прежнему конечен. 
Для $z\in \CC$  при $r\in \RR^+$ положим 
\begin{equation}\label{B}
{\sf B}_u(z,r):=\frac{1}{\pi r^2}\int_{|z'-z|\leq r} u(z')\dd \lambda (z'), \quad 
{\sf B}_u^{\text{rad}}(r):={\sf B}_u(z,r) 
\end{equation}
 --- соответственно средние  функции $u$ по замкнутым кругам  $\overline  D(z,r):=\{z' \in \CC \colon |z'-z|\leq r\}$ и $\overline  D(r):=\overline  D(0,r)$,  а ${\sf M}_u(S):=\sup_S u$ --- верхняя грань $u$ на  $S\subset \CC$. Для функции $u$ конечного порядка функция ${\sf B}_u^{\text{rad}}$ также функция конечного порядка в том смысле, что  после замены  $\ln|f|(z)$ на ${\sf B}_u^{\text{rad}}\bigl(|z|\bigr)$ верхний предел в \eqref{fo} по-прежнему конечен.

По условию ограниченности сверху   на $\CC\!\setminus\!E$ функции $u$ 
можно рассмотреть субгармоническую функцию  $u-{\sf M}_u\bigl(\CC\!\setminus\!E\bigr)$, отрицательную на $\overline \CC\!\setminus\!E$, сохранив для неё то же обозначение $u$.  Тогда 
\begin{multline}\label{E}
{\sf B}_u^{\rm rad}(r)\leq {\sf B}_{u^+}^{\rm rad}(r)\leq \frac{1}{\pi r^2}\biggl(\int_{\overline D(r)\!\setminus\!E}+
\int_{\overline D(r)\cap E}\biggr)u^+\dd \lambda\\
\leq \frac{1}{\pi r^2}\int_{E\cap \overline D(r)}{\sf B}_{u^+}(z,r)\dd \lambda (z)
\end{multline}
ввиду неравенств $u^+(z)\leq {\sf B}_{u^+}(z,r)$ при всех $z\in \CC$ для субгармонической функции $u^+$. 
Но при $z\in \overline D(r)$ из включений $\overline D(z,r)\subset \overline D(2r)$ для {\it положительной\/} функции $u^+$ имеем
\begin{equation*}
{\sf B}_{u^+}(z,r)= \frac{1}{\pi r^2} 
\int_{\overline D(z,r)}u\dd \lambda
\overset{\eqref{B}}{\leq} \frac{1}{\pi r^2} 
\int_{\overline D(2r)}u\dd \lambda=4{\sf B}_{u^+}^{\rad}(2r)\text{ при $z\in \overline D(r)$}.
\end{equation*}
Отсюда, продолжая \eqref{E}, получаем 
\begin{multline*}
{\sf B}_{u^+}^{\rm rad}(r)\leq \frac{1}{\pi r^2}\int_{E\cap \overline D(r)}4{\sf B}_{u^+}^{\rad}(2r)\dd \lambda (z)\\
\leq 4{\sf B}_{u^+}^{\rad}(2r) \frac{\lambda (E\cap \overline D(r))}{\pi r^2}
\overset{\eqref{m00}}{=}o(1){\sf B}_{u^+}^{\rad}(2r)\text{ при $r\to +\infty$}.
\end{multline*}
 Возрастающая функция ${\sf B}_{u^+}^{\rm rad}$  конечного порядка  на положительной полуоси при таком свойстве   --- тождественный  нуль и, как следствие, $u^+\equiv 0$ на $\CC$. Вспоминая переобозначение для $u$, получаем $u\equiv{\sf M}_u(\CC \!\setminus\!E)$.  
\end{proof}  
\begin{remark}
Доказательство  даёт теорему Лиувилля и  для субгармонических функций $u$ конечного порядка, ограниченных сверху на $\CC\!\setminus\!E$.
\end{remark}
 Приведённое доказательство имеет значительный потенциал развития и обобщения в нескольких  направлениях: 
\begin{itemize}
\item для целых, субгармонических и плюрисубгармонических функций на плоскости и в пространстве, а также для выпуклых функций, начиная с функций на прямой;
\item вместо усреднений  по кругам из  \eqref{B} можно рассмотреть усреднения по окружностям, сферам, по гармоническим мерам на замкнутых контурах, уходящих в $\infty$, и, максимально общ\'о, по мерам Йенсена в сочетании с развитием неравенств типа Харнака для {\it положительных\/} субгармонических функций и с привлечением утверждений типа лемм о малых дугах Эдрея\,--\,Фукса;
\item для  голоморфных или (плюри)субгармонических функций в круге или шаре при условии достаточно быстрого стремления к нулю меры пересечения  исключительного множества $E$ с кольцами или шаровыми слоями, сужающимися при приближении  к границе, но с заключением об ограниченности функции или ограничении на её рост, а не о постоянстве;   
\item с учётом баланса между малостью исключительного множества $E$ и интегралами от  функции на дополнении исключительного множества $E$, а также глобальным ростом функции на всей области определения.
\item ограничения можно рассматривать не по всем кругам, шарам, сферам, контурам и т.п., а по их достаточно редким последовательностям в балансе с ограничениями на рост функций и малость $E$.  
  
\end{itemize}
  
Автор  глубоко признателен А.\,Д.~Баранову, благодаря чьему пленарному докладу  и стимулирующим on-line контактам с ним на Международной научной конференции <<Комплексный анализ и его приложения>> в г. Казань (24--28 августа 2020 г.)  появились как эта заметка, так и перспективные планы по развитию её содержания.

{\sl  Хабибуллин Булат Нурмиевич}

 {Башкирский государственный университет, г. Уфа}

{E-mail: {\tt khabib-bulat@mail.ru}}
\end{document}